\documentclass[a4paper,twoside,12pt]{article} 
\usepackage{amsfonts,amssymb,mathrsfs,amsmath}
\usepackage[dvips]{color}
\usepackage{graphicx}
\title{On the structure of cube tiling codes}
\date{}
\author{Andrzej P. Kisielewicz\\
\\
{\small Wydzia{\l} Matematyki, Informatyki i Ekonometrii, Uniwersytet Zielonog\'orski}\\
{\small ul. Z. Szafrana 4a, 65-516 Zielona G\'ora, Poland}\\
{\small A.Kisielewicz@wmie.uz.zgora.pl}\\
}
\numberwithin{equation}{section}
\newtheorem{pr}{\sc Proposition}

\newtheorem{lemat}[pr]{\sc Lemma}

\newtheorem{tw}[pr]{\sc Theorem}

\newtheorem{df}{\sc Definition}

\newtheorem{uw}{\sc Remark}
\newtheorem{uwi}[uw]{\sc Remarks}

\newtheorem{nap}{\sc Example }

\newtheorem{nps}[nap]{\sc Examples}

\def\ka #1{\mathscr{#1}}
\def\kal #1 #2{\mathscr{#1}^{#2}}
\def\proof{\noindent \textit{Proof.\,\,\,}}

\def\zet{\mathbb{Z}}

\def\er{\mathbb{R}}
\def\en{\mathbb{N}}

\def\te{\mathbb{T}}

\def\iver #1{\mbox{\tt [} #1 \mbox{\tt]}}

\begin{document}

\numberwithin{pr}{section}
\numberwithin{uw}{section}
\maketitle
\begin{abstract}

Let $S$ be a set of arbitrary objects, and let $S^d=\{v_1...v_d\colon v_i\in S\}$.  A polybox code is a set $V\subset S^d$ with the property that for every two words $v,w\in V$ there is $i\in [d]$ with $v_i'=w_i$, where a permutation $s\mapsto s'$ of $S$ is such that $s''=(s')'=s$ and $s'\neq  s$. 
If $|V|=2^d$, then $V$ is called a cube tiling code. Cube tiling codes determine  $2$-periodic cube tilings of $\er^d$ or, equivalently, tilings of the flat torus $\te^d=\{(x_1,\ldots ,x_d)({\rm mod} 2):(x_1,\ldots ,x_d)\in \er^d\}$ by translates of the unit cube as well as  $r$-perfect codes in $\zet^d_{4r+2}$ in the maximum metric. By a structural result, cube tiling codes for $d=4$ are  enumerated. It is computed that there are 27,385 non-isomorphic cube tiling codes in dimension four, and the total number of such codes is equal to $17,794,836,080,455,680$. Moreover, some procedure of passing from a cube tiling code to a cube tiling code in dimensions $d\leq 5$ is given.

\medskip
\noindent
\textit{Key words:} box, code, cube tiling, perfect code.

\end{abstract}

\section{Introduction}
Let $S$ be a set of arbitrary objects which  will be called an \textit{alphabet}, and the elements of $S$ will be called \textit{letters}. A permutation $s\mapsto s'$ of the alphabet $S$ such that $s''=(s')'=s$ and $s'\neq  s$ is said to be a \textit{complementation}. Through the paper we use a fixed complementation and therefore the alphabet $S$ will be given in the form $S=\{a_1,a_1',...,a_k,a_k'\}$. 
Let $S^d=\{v_1...v_d\colon v_i\in S\}$. Elements of the set $S^d$ are called {\it words}. A {\it polybox code}  (or  a {\it polybox genom}) is a set $V\subset S^d$ with the property that for every two words $v,w\in V$ there is $i\in [d]=\{1,...,d\}$ with $v_i'=w_i$ (\cite{KP}). If $|V|=2^d$, then $V$ will be called a {\it cube tiling code} (\cite{LS2}). 

  


\smallskip
A natural model for a polybox code is some special system of boxes, where a {\it box } is a set $A\subseteq X=X_1\times\cdots \times X_d$ of the form $A=A_1\times\cdots \times A_d$, where $A_i\subseteq X_i$ for each $i\in [d]$. We assume that $|X_i|>1$ for $i \in [d]$ and  call the box $X$ a {\it $d$-box.} The mentioned system of boxes can be derived from a code in the following way: 
For $i\in [d]$ let $f_i\colon S\to 2^{X_i}\setminus \{\emptyset,X_i\}$ be such that $f_i(s')=X_i\setminus f_i(s)$, and let $f\colon S^d\to 2^{X}$ be defined by the formula $f(v_1\ldots v_d)=f_1(v_1)\times\cdots\times f_d(v_d)$. If now $V\subseteq S^d$ is a polybox code, then the set $f(V)=\{f(v)\colon v\in V\}$ is a system of pairwise disjoint boxes in $X$ such that for every two boxes $K,G\in f(V)$ there is $i\in [d]$ with $K_i=X_i\setminus G_i$ (compare Figure 1). We call the set $f(V)$ a {\it realization} of $V$. Boxes $K,G$ with this last property are called {\it dichotomous}. Similarly, two words $v,w\in S^d$ are {\it dichotomous} if $v_i=w_i'$ for some $i\in [d]$. About the above defined function $f$ we say that it {\it preserves dichotomies}.

A $2$-periodic cube tiling of $\er^d$ can be treated as a realization of a cube tiling code. Recall that a family  $[0,1)^d+T=\{[0,1)^d+t\colon t\in T\}$ is a {\it cube tiling} of $\er^d$ if  every two cubes in $[0,1)^d+T$ are disjoint and $\bigcup_{t\in T} ([0,1)^d+t)=\er^d$. A tiling $[0,1)^d+T$ is called {\it $2$-periodic} if $T+2\zet^d=T$. To obtain a $2$-periodic cube tiling from a cube tiling code $V\subset S^d$ let $P\subset [0,1)$ be a set with $k$ elements, and let $f_i(a_j)=[0,1)+p_j+2\zet, f_i(a_j')=[1,2)+p_j+2\zet$, where $p_j\in P$ for $j\in [k]$, and $i\in [d]$.   
The realization $f(V)$ of $V$ is a partition of the $d$-box $\er^d$ into $2^d$ pairwise dichotomous boxes $f(v)$, $v\in V$. It is easy to see that at the same time it can be viewed as a $2$-periodic cube tiling of $\er^d$ (Figure 1). Clearly, the tiling $f(V)$ determines a tiling of the flat torus $\te^d=\{(x_1,\ldots ,x_d)({\rm mod} 2):(x_1,\ldots ,x_d)\in \er^d\}$ by translates of the unit cube as well as an $r$-perfect code in $\zet^d_{4r+2}$ in the maximum metric (Figure 1). Recall that, a set $C\subset \zet^d_n$ is an {\it $r$-perfect code $($in the metric $\delta$$)$} if for every $x\in \zet_n^d$ there is exactly one $c\in C$ such that $\delta(c,x)\leq r$. In other words, the family  $\{B(c,r)\colon c\in C\}$ consists of mutually disjoint balls $B(c,r)=\{x \in \zet^d_n\colon \delta(c,x)\leq r\}$ and $\zet^d_n=\bigcup_{c\in C} B(c,r)$. Every cube tiling code $V\subseteq S^d$ induces an $r$-perfect code $C\subset \zet^d_{4r+2}$, $r\in \en$, in the {\it maximum} metric $\delta_\infty(x,y)=\max_{1\leq i\leq d}|x_i-y_i|$ (compare \cite{Co}). Indeed, consider a realization $f(V)$, where $f:S^d\to 2^{\zet^d_{4r+2}}$ preserves dichotomies and for every $i\in [d]$ and every $s\in S$ the set $f_i(s)$ consists of $2r+1$ consecutive integers in $\zet_{4r+2}$. Obviously, the set of pairwise dichotomous boxes $f(V)$ is a partition of $\zet^d_{4r+2}$. Since each box in $f(V)$ is a ball of radius $r$ in the metric $\delta_\infty$, taking the centers of each box from $f(V)$ we obtain an $r$-perfect code $C$. 

\vspace{-0mm}
{\center
\includegraphics[width=13cm]{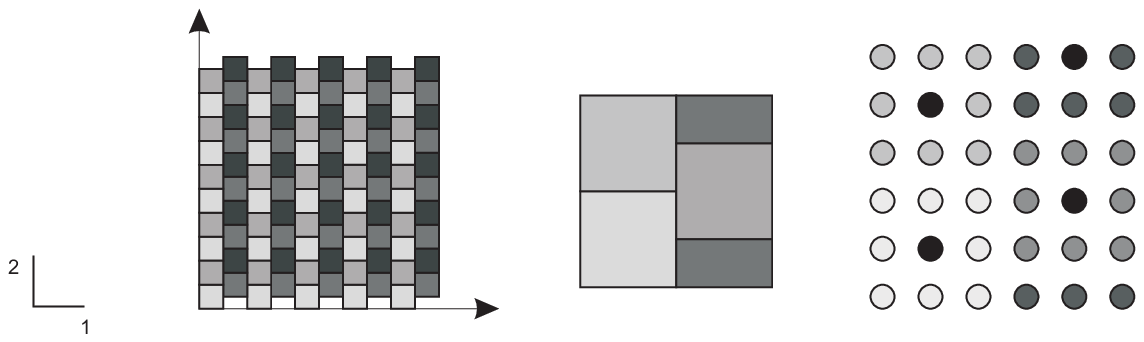}\\
}

\medskip
\noindent{\footnotesize Figure 1: Three realizations $f^i(V)$, $i=1,2,3$, of the cube tiling code $V=\{aa,aa',a'b,a'b'\}$. The realization $f^1(V)$, where $f_1^1(a)=[0,1)+2\zet, f_2^1(a)=[0,1)+2\zet$ and $f_2^1(b)=[0,1)+\frac{1}{2}+2\zet$, is a $2$-periodic cube tiling of $\er^2$ (picture on the left); the realization $f^2(V)$, where $f_1^2(a)=[0,1), f_2^2(a)=[0,1)$ and $f_2^2(b)=[\frac{1}{2},\frac{3}{2})$ is a cube tiling of $\te^2$ (picture in the middle). Finally, the realization $f^3(V)$, where $f_1^3(a)=\{0,1,2\}, f_2^3(a)=\{0,1,2\}$ and $f_2^3(b)=\{1,2,3\}$, is a tiling of $\zet^2_6$ by boxes (balls in the maximum metric), (picture on the right). Corresponding $1$-perfect code in the maximum metric is represented by the black dots. 
}

\medskip

In the presented paper we enumerate all cube tiling codes $V\subset \{a_1,a_1',...,a_8,a_8'\}^4$ which leads to a complete classification of cube tiling codes in dimension four and thus, a complete classification of $2$-periodic cube tilings of $\er^4$, cube tilings of $\te^4$ as well as  $r$-perfect codes in $\zet^d_{4r+2}$ in the maximum metric.   Complete classifications of cube tiling codes $V\subset \{a,a',b,b'\}^d$ for $d\leq 5$ are, so far, known: For $d\leq 4$ it was computed by Dutour and Itoh in \cite{DIP}, and for $d=5$ by Mathew, {\"O}sterg{\aa}rd and Popa in \cite{Ost}. To classify cube tiling codes $V\subset S^4$, where $S$ is an arbitrary alphabet, we show that every cube tiling code in dimension four can be merged into one of twenty one codes of some special form (Theorem \ref{kod}). This in turn reveals some meta-structures of $2$-periodic cube tilings of $\er^4$ in the spirit of \cite{LP2}. 

One of the most known structural problem dealing with cube tilings of $\er^d$ is Keller's conjecture (\cite{CS2, De,Ke1,KisL, Kap, Kis,LS1,LP,M,P}). This conjecture, which says that in every cube tiling of $\er^d$ there is a pair of cubes having a complete facet in common, is however only one of many interesting problems on the structure of cube tilings (\cite{DI,DIP,DI2,GKP,H,KP,LS2,SS,Sz}).    
At the end of the paper we examine one of such structural problem which was previously considered in \cite{DI2,Ost} for cube tiling codes  $V\subset \{a,a',b,b'\}^d$ for $d\leq 5$ (in the language of cube tilings of $\te^d$). It concerns a procedure of passing from a code to a code via a sequence of local transformations of codes. In the paper we generalize such procedure on arbitrary cube tiling codes $V\subset S^d$ in dimensions up to five (Theorem \ref{cut}).


\medskip
A non-empty set $F\subseteq X$ is said to be a \textit{ polybox} if
there is a set of pairwise dichotomous boxes $\ka F$ such that $\bigcup \ka F=F$. The set $\ka F$ is called a {\it suit} for $F$. 
We add an extra letter $\ast$ to the set $S$ and the set $S\cup \{*\}$ is denoted by $*S$. We assume that $*'=*$ and the star is the only letter with this property. Two words $v,w\in (*S)^d$ are {\it dichotomous} if $v_i=w_i'$ for some $i\in [d]$, where $v_i,w_i\in S$, and  $V\subset (*S)^d$ is a  \textit{polybox code} if it consists of pairwise dichotomous words. A word containing $*$ is called {\it improper}. The meaning of $*$ is that $f_i(*)=X_i$,  while $f_i(s)\in 2^{X_i}\setminus \{X_i,\emptyset\}$ for every $s\in S$ and every $i\in [d]$.
A set $V\subset (*S)^d$ of pairwise dichotomous words is said to be a {\it partition code} if any realization $f(V)$ of $V$ is a suit for a $d$-box $X$. 
Codes $V,W\subset (*S)^d$ are {\it equivalent}, which is denoted by $V\equiv U$,  if $\bigcup f(V)=\bigcup f(W)$ for every $f$ that preserves dichotomies (see Figure 4 and 5). 

If $v\in (*S)^d$, and $\sigma$ is a permutation of the set $[d]$, then $\bar{\sigma}(v)=v_{\sigma(1)}\ldots v_{\sigma(d)}$. For every $i\in [d]$ let $h_i:*S\rightarrow *S$ be a bijection such that $h_i(l')=(h_i(l))'$ for every $l\in *S$ and $h_i(*)=*$. Such $h_i$ will be called a {\it position bijection} ({\it at position i}).   Let $h:(*S)^d\rightarrow (*S)^d$ be defined by the formula $h(v)=h_1(v_1)\ldots h_d(v_d)$. 
The group of all possible mappings $h\circ \bar{\sigma}$  will be denoted by $G((*S)^d)$ or $G(S^d)$ depending on whether we consider words written down in the alphabet $*S$ or in $S$. Let $S$ and $T$ be two alphabets with  complementations, and let $|S|\leq |T|$. Two polybox codes $V\subset (*S)^d$ and $U\subset (*T)^d$ are {\it isomorphic} if there is  $h\circ \bar{\sigma}\in G((*T)^d)$ such that $U=h\circ \bar{\sigma}(i(V))$, where $i\colon *S\rightarrow *T$ is a fixed injection such that $i(*)=*$ and $i(s')=i(s)'$ for $s \in S$. The composition $h\circ \bar{\sigma}$ is an {\it isomorphism} between $V$ and $U$. 
By a set of {\it all non-isomorphic cube tiling codes in dimension $d$} we mean any set $\ka N_d\subset S^d$ such that for every cube tiling code $W\in P^d$ there is $V\in \ka N_d$ such that $V$ and $W$ are isomorphic. 
 

\section{Enumeration of twin pair free partition codes}

Two words $v,u \in (*S)^d$ are called  a {\it twin pair} if there is $i \in [d]$ such that $v_i=u_i'$, $v_i\neq *$, and $v_j=u_j$ for every $j\in [d]\setminus \{i\}$. We say that a twin pair $v,u$ is {\it glued} at the $i$th position if the pair $v,u$ is replaced by the improper word $w$ having  the star at the $i$th position, where $v_i=u_i'$ and $w_j=v_j$ for all $j\neq i$. The word $w$ is called a  {\it gluing} of words $v,u$.  
To every polybox code $V$ we can assign its twin pair free code $F_V$ which arises from $V$ by successive gluing of words which form a twin pair. (Such assignment is not usually unique.)  More precisely, in the first step we glue a twin pair in $V$, say it is $v_1v_2v_3v_4$ and $v_1v_2v_3v_4'$, obtaining the code $V^1=V\setminus \{v_1v_2v_3v_4,v_1v_2v_3v_4'\}\cup \{v_1v_2v_3*\}$ (if there is no twin pair in $V$ we take $F_V=V$). Note that $V^1$ is a polybox code. If $V^1$ does not contain a twin pair we take $F_V=V^1$. If there is a twin pair in $V^1$ we proceed as above obtaining a code $V^2$. After $n$ steps, where $n\leq 2^{d}-1$ and $V^i$ contains a twin pair for $i\in [n-1]$, we obtain a twin pair free code $V^n$ and then $F_V= V^n$.  For example, if $V=\{aaa, \;aba',\;ab'a',   \;a'ab,\;a'ab',   \;a'a'a', \;ba'a, \;b'a'a\}$, then  $V^1=\{aaa, \;a*a',  \;a'ab,\;a'ab',   \;a'a'a', \;ba'a, \;b'a'a\}$, $V^2=\{aaa, \;a*a',  \;a'a*,   \;a'a'a', \;ba'a, \;b'a'a\}$, $V^3=\{aaa, \;a*a',  \;a'a*,   \;a'a'a', \;*a'a\}$, and then $F_V=V^3$.

\vspace{-0mm}
{\center
\includegraphics[width=12cm]{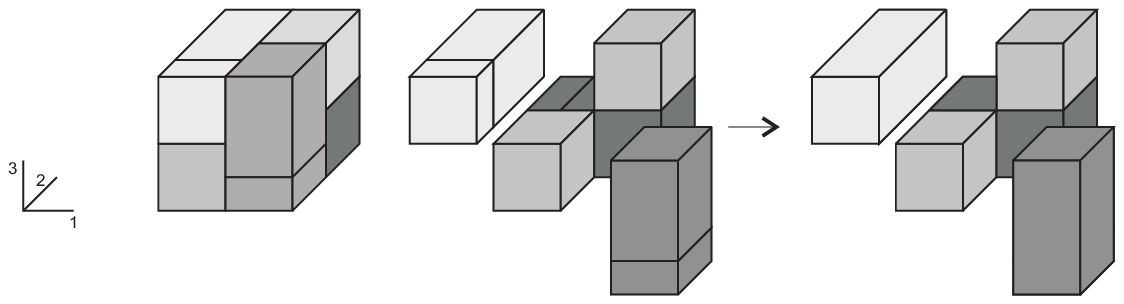}\\
}

\medskip
\noindent{\footnotesize Figure 2: A realization of the cube tiling code $V=\{aaa, \;aba',\;ab'a',   \;a'ab,\;a'ab',   \;a'a'a', \;ba'a, \;b'a'a\}$ (on the left) and a realization of the code $F_V=\{aaa, \;a*a',  \;a'a*,   \;a'a'a', \;*a'a\}$ (on the right).
}

\medskip

It can happen that $F_V=\{****\}$. This however means that, up to isomorphism, $V^{n-1}=\{v_1***,v_1'***\}$. This, in turn, means that the code $V$ (as well as $V^{n-1}$) belongs to the class of layered codes: A code $U\subset (*S)^d$ is {\it layered} if there are $s\in S$ and  $i\in [d]$ such that $u_i\in \{s,s'\}$ for every word $u\in U$. 

\bigskip
To enumerate all twin pair free partition codes $V\subset (*S)^4$ we shall use an algorithm which is described in \cite{Kis}. 
Let $v\in (*S)^d$, and let
\begin{equation}
\label{nor}
|v|=\prod^d_{i=1}(2\iver {v_i=*}+\iver {v_i\neq*}),
\end{equation}
where $\iver {p}=1$ if the sentence $p$ is true and $\iver {p}=0$ if it is false. It can be proved that a polybox code $V\subset (*S)^d$ is a partition code if and only if $\sum_{v\in V}|v|=2^d$ (\cite{KP}). 
Thus, denoting by $x_i$ the number of words $v$ in  a partition code $V$ such that $|v|=2^i$, $i\in \{0,1,2,3\}$, to every code $V$ we may assign the system of equations  $\sum_{i=0}^{3}x_i2^i=2^4$, $\sum_{i=0}^{3}x_i=k$, where $k$ is the number of words in $V$. Moreover, in \cite{KP} we showed that every partition code $V\subset (*S)^d$ having at least two words contains a pair of words $v,u$, such that $v_i=u_i$ or $v_i=u_i'$ for every $i\in [d]$ and 
\begin{equation}
\label{pa}
|\{i\colon v_i=u_i'\;\;{\rm and}\;\; v_i\neq *\}|\equiv 1\;\; ({\rm mod}\; 2).
\end{equation}

For the reasons that will be explained below, it is enough to consider the alphabet $*S=\{a,a',b,b',*\}$.

\medskip
{\bf Algorithm 1}.

\medskip
Let $*S=\{a,a',b,b',*\}$, $V_{3,0}=\{aaaa,a'a'a'a\}$ and $V_{3,1}=\{aaa*,a'a'a'*\}$. By (\ref{pa}), every twin pair free partition code $V\subset (*S)^4$ contains, up to isomorphism, the code $V_{3,0}$ or $V_{3,1}$. Let $\ka C^k$ be the family of all $k$-elements twin pair free partition codes $V\subset (*S)^4$ which contain the code $V_{3,0}$ or $V_{3,1}$.  Our goal is to find the family $\ka C^k$.

\smallskip
{\it Input}. The number $k$.

{\it Output}. The family $\ka C^k$.

\smallskip
1. Let $\ka S_k=\{(x_0,...,x_{3}) \in \en^4\colon  \sum_{i=0}^{3}x_i2^i=2^4\; {\rm and}\; \sum_{i=0}^{3}x_i=k\}$, where $\en=\{0,1,2,...\}$.

2. For $i \in \{0,...,3\}$ indicate the set $\ka A_i$ consisting of all words $v\in (*S)^4$ such that $v$ contains precisely $i$ letters $*$. 


3. Fix $x\in \ka S_k$ and let $s(x)=\{i_1<\cdots <i_m\}$ consists of all $i_j\in \{0,...,3\}, j\in [m]$, for which $x_{i_j}>0$. Fix $V_{3,i_1}\subset \ka A_{i_1}$. 
For $i\in s(x)$ let $\ka B_i=\{v\in \ka A_i\colon V_{3,i_1}\cup \{v\}\;\; {\rm is\; a \; twin\;\; pair\;\; free\;\; code}\}$.


4. Let $I$ be the multiset containing $i_1$ with the multiplicity $x_{i_1}-2$ (recall that $x_{i_1}\geq 2$) and  $i_j$ with the multiplicity $x_{i_j}$ for $j\in \{2,...,m\}$. By $I[j]$ we denote the $j$th element of $I$. 

5. Let $\ka D^2=\{V_{3,i_1}\}$. 

6. For $l\in \{2,...,k-1\}$ having computed $\ka D^l$ we compute the set $\ka D^{l+1}$: For $v\in \ka B_{I[l]}$ and for $U\in \ka D^l$ if $U\cup \{v\}$ is a twin pair free code, then we attach it to  $\ka D^{l+1}$. 

7. Clearly, $\ka D^k=\ka D^k(V_{3,i_1},x)$ so let $\ka C^k$ be the union of the sets $\ka D^k(V_{3,0},x)\cup \ka D^k(V_{3,1},x)$ over $x\in \ka S^k$. 



\medskip
To explain why in the above algorithm we may assume that $*S=\{a,a',b,b',*\}$, which simplifies computations, we need to take a look at the structure of a partition code $V\subset (*S)^d$. If $l\in *S$ and $i\in [d]$, then let $V^{i,l}=\{v\in V\colon v_i=l\}$. If $S=\{a_1,a'_1,...,a_k,a_k',*\}$, then the representation 
$$
V=V^{i,a_1}\cup V^{i,a_1'}\cup \cdots \cup V^{i,a_{k}}\cup V^{i,a_{k}'}\cup V^{i,*}
$$
will be called a {\it distribution of words } in $V$. For example, if $V=\{aaaa, \;aaa'b, \;aa'a'*, \;a*aa', \\a'aab', \;a'a*b, \;a'a'a'a, \;a'a'*a', \;*aa'b', \;*a'aa\}$, then 
$$
V^{4,a}=\{aaaa, \;a'a'a'a,\;*a'aa\},\;\; V^{4,a'}=\{\;a*aa',\;a'a'*a'\},
$$
$$
V^{4,b}=\{aaa'b, \;a'a*b\}, \;\;V^{4,b'}=\{a'aab',\;*aa'b'\}, \;\; V^{4,*}=\{aa'a'*\}.
$$
Any realization of the set $V^{i,l}\cup V^{i,l'}$, where $V$ is a partition code, is a {\it cylinder} (Figure 2). It is a set $F\subseteq X$ such that  for every set $l_i=\{x_1\}\times \cdots \times \{x_{i-1}\}\times X_i \times \{x_{i+1}\}\times \cdots \times \{x_d\}$, where $x_j\in X_j$ for $j\in [d]\setminus\{i\}$, 
one has 

$$
l_i\cap F=l_i \; \; \; {\rm or}\; \; \; l_i\cap F=\emptyset,
$$

\vspace{-0mm}
{\center
\includegraphics[width=9cm]{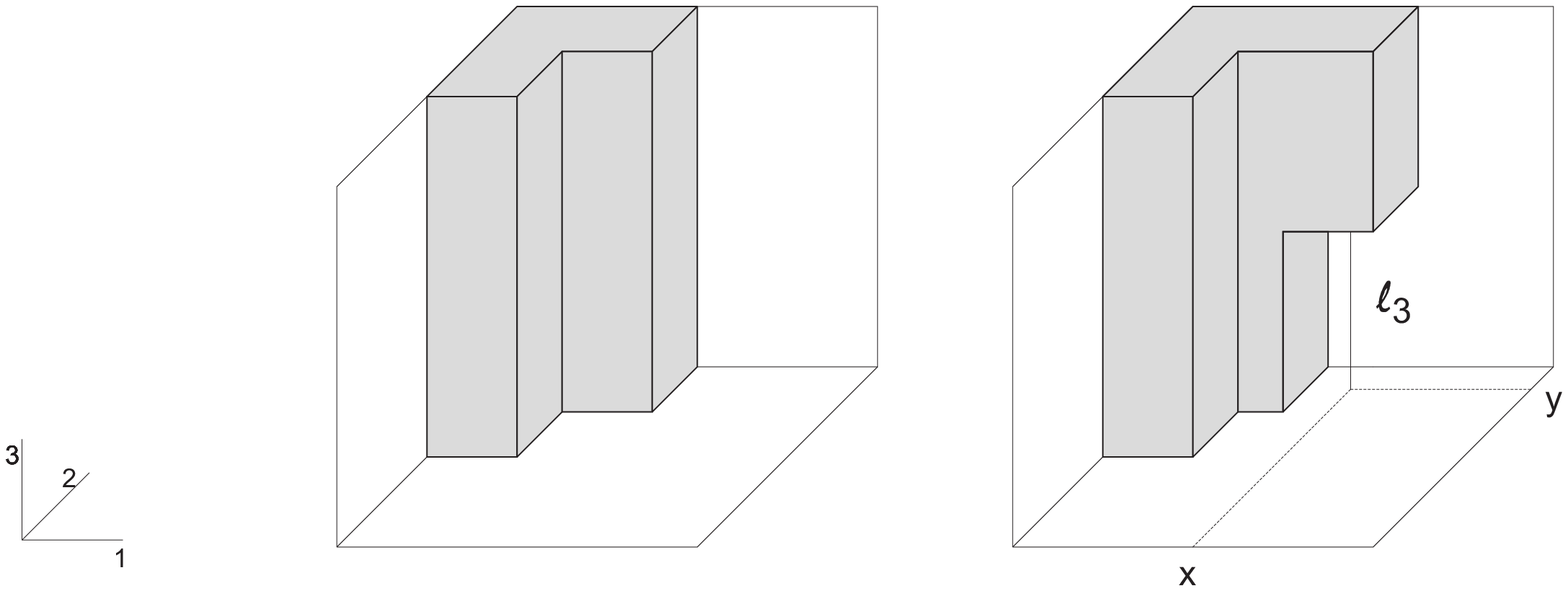}\\
}

\medskip
\noindent{\footnotesize Figure 3: The set on the left is a $3$-cylinder in $X=[0,1]^3$, and the set on the right is not, because the set $l_3=\{x\}\times \{y\}\times [0,1]$ has a non-empty intersection with this set but $l_3$ is not entirely contained in it. 
}

\medskip
Since any realization of $V^{i,l}\cup V^{i,l'}$, where $V$ is a partition code, is a cylinder we shall say that a partition code, in particular a cube tiling code, has the {\it cylindrical structure}.

Observe that, for every cube tiling code $V\subset S^d$ we have $|V^{i,a_{k}}\cup V^{i,a_{k}'}|\geq 2$, whenever $V^{i,a_{k}}\cup V^{i,a_{k}'}\neq\emptyset$. Thus, we get

\begin{lemat}
\label{uu}
To obtain a set $\ka N_d\subset S^d$ of all non-isomorphic cube tiling codes in dimension $d$ it is enough to take $S=\{a_1,a_1',...,a_{2^{d-1}},a_{2^{d-1}}'\}$.  
\hfill{$\square$}
\end{lemat}

To justify the assumption $S=\{a,a',b,b',*\}$ let $V^{i,l}_{i^c}=\{v_1...v_{i-1}v_{i+1}...v_d\colon v\in V^{i,l}\}$.
Since every realization of $V^{i,l}\cup V^{i,l'}$ is a cylinder, it follows that  realizations of the sets $V^{i,l}_{i^c}$ and $V^{i,l'}_{i^c}$ are equal, which means, that the codes $V^{i,l}_{i^c}, V^{i,l'}_{i^c}\subset (*S)^{d-1}$ are equivalent. It is easy to check that if two disjoint codes $V,U\subset (*S)^3$ are twin pair free and equivalent then, there there are improper words in both codes. Moreover, if  a code $V^{i,l}\cup V^{i,l'}\subset (*S)^4$ does not contain a twin pair, then $|V^{i,l}\cup V^{i,l'}|\geq 4$. Therefore, if $V\subset (*S)^4$ is a twin pair free partition code and  $V^{i,l}\cup V^{i,l'}\neq \emptyset$ for $l\in \{a,b,c\}$, then $|V|\geq 12$. Since every solution of the system of equations $\sum_{i=0}^{3}x_i2^i=2^4$ and $\sum_{i=0}^{3}x_i=k$ for $k\geq 12$ has the property that $\sum_{i=1}^{3}x_i\leq 4$, it follows that in every partition code  $V\subset (*S)^4$ with $|V|\geq 12$ there are at most four improper words. On the other hand, as it was noted above, for every $l\in \{a,a',b,b',c,c'\}$ the set $V^{i,l}$ contains an improper word. Thus, there are at least six improper words in $V$, which is a contradiction.  

\begin{tw}
\label{kod}

There are, up to isomorphism, twenty twin pair free partitions codes in dimension four. These are:
  
\begin{enumerate}
\item $C^1=\{a*aa, \;a**a', \;a'*a*, \;a'*a'a', \;**a'a\},$
\item $C^2=\{aaa'a, \;aa*a', \;a'aaa', \;a'aa'*, \;*aaa, \;*a'**\},$
\item $C^3=\{a*aa', \;a**a, \;a'aaa', \;a'a*a, \;a'a'a*, \;a'a'a'a, \;**a'a'\},$
\item $C^4=\{aaaa', \;a*a'*, \;a'aa'a, \;a'a*a', \;a'a'a'*, \;*a'aa', \;**aa\},$
\item $C^5=\{aa*a', \;aa'a'*, \;a*aa, \;a'aa*, \;a'a'*a, \;a'*a'a', \;*aa'a, \;*a'aa'\},$
\item $C^6=\{aaaa, \;aaa'*, \;aa'*a', \;a'a*a, \;a'a'aa', \;a'*a'a', \;*aaa', \;*a'*a\}$
\item $C^7=\{aaaa, \;aa'a'a', \;aa'*a, \;a*aa', \;a'aaa', \;a'a'a'a, \;a'a'*a', \;a'*aa, \;*aa'*\},$
\item $C^8=\{aaaa, \;aa'*a, \;a'a'a'a, \;a'*aa, \;*aa'*, \;ba'a'a', \;b*aa', \;b'aaa', \;b'a'*a'\},$
\item $C^9=\{aaaa, \;aa'*a, \;a**a', \;a'a'a'a, \;a'baa, \;a'b*a', \;a'b'a* a'b'a'a', \;*aa'a\},$
\item $C^{10}=\{aaaa, \; a*a'a, \;a'a*a, \;a'a'a'a, \;ba*a', \;ba'a*, \;ba'a'a', \;b'a'aa, \;b'**a'\}$
\item $C^{11}=\{aaaa, \; aa*a', \;aa'a'b, \;aa'*b', \;a'aab, \;a'a'a'a, \;a'*ab', \;a'*a'a', \;*aa'a, \;*a'ab\},$
\item $C^{12}=\{aaaa, \; a*a'a, \;a'a*a, \;a'a'a'a, \;*a'aa, \;*baa', \;b*a'a', \;bb'aa', \;b'ba'a', \;b'b'*a'\},$
\item $C^{13}=\{aaaa, \; a*a'a, \;a'a*a, \;a'a'a'a, \;*a'aa, \;*a'a'a', \;ba*a', \;ba'aa', \;b'aa'a', \;b'*aa'\},$
\item $C^{14}=\{aaaa, \; a*a'a, \;ab*a', \;ab'b'a', \;a'a*a, \;a'a'a'a, \;a'*b'a', \;a'bba', \;*a'aa, \;*b'ba'\},$
\item $C^{15}=\{aaaa, \; a*a'a, \;a'a*a, \;a'a'a'a, \;*a'aa, \;*bb'a', \;b*ba', \;bb'b'a', \;b'bba', \;b'b'*a'\},$
\item $C^{16}=\{aaaa, \;aaa'b, \;aa'a'*, \;a*aa', \;a'aab', \;a'a*b, \;a'a'a'a, \;a'a'*a', \;*aa'b', \;*a'aa\},$
\item $C^{17}=\{aaaa, \;aa'aa', \;aa'*a, \;a*a'a', \;a'aa'a', \;a'a'a'a, \;a'a'*a', \;a'*aa, \;*aaa', \;*aa'a\},$
\item $C^{18}=\{aaaa, \; a*a'a, \; abaa', \;ab'*a', \;a'a*a, \;a'a'a'a, \;a'*aa', \;a'b'a'a', \;*a'aa, \;*ba'a'\},$
\item $C^{19}=\{aaaa, \; aab'a', \; aa'ba, \;aa'b'*, \;abba', \;a'a'a'a, \;a'*aa, \;a'b*a', \;a'b'b'a', \;*aa'a, \;*b'ba'\},$
\item $C^{20}=\{aaaa, \;aaba', \; aa'b*, \;aa'b'a, \;ab'b'a', \;a'a'a'a, \;a'baa, \;a'bba', \;a'b'a*, \;a'b'a'a', \;*aa'a, \;*bb'a'\}.$
\end{enumerate}
\hfill{$\square$}
\end{tw}

Each of these twenty codes can be partially visualized. 
For example, let $C=C^{16}$. We have
$$
C^{4,a}=\{aaaa, \;a'a'a'a,\;*a'aa\},\;\; C^{4,a'}=\{\;a*aa',\;a'a'*a'\},
$$
$$
C^{4,b}=\{aaa'b, \;a'a*b\}, \;\;C^{4,b'}=\{a'aab',\;*aa'b'\}, \;\; C^{4,*}=\{aa'a'*\}.
$$

The two codes $C^{4,a}_{4^c}=\{aaa, \;a'a'a',\;*a'a\}$ and   $C^{4,a'}_{4^c}=\{\;a*a,\;a'a'*\}$ are equivalent and similarly the codes $C^{4,b}_{4^c}=\{aaa', \;a'a*\}$ and $C^{4,b'}_{4^c}=\{a'aa,\;*aa'\}$ are equivalent. Below we draw a realizations of $C$. 

\vspace{-0mm}
{\center
\includegraphics[width=11cm]{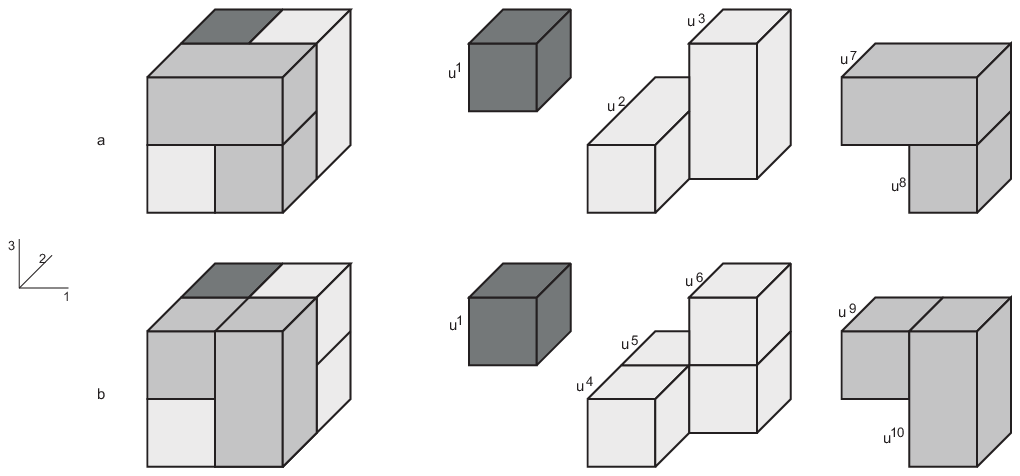}\\
}

\medskip
\noindent{\footnotesize Figure 4: In the figure ($a$) we see a realization of the code $(C^{4,a}\cup C^{4,b}\cup C^{4,*})_{4^c}$ (that is, the code $C^{4,a}\cup C^{4,b}\cup C^{4,*}$ with the fourth coordinate removed) and similarly, in ($b$) we have a realization of the code $(C^{4,a'}\cup C^{4,b'}\cup C^{4,*})_{4^c}$. The realization $f(C)$ of $C=C^{16}$ is such that $\bigcup f(C)=[0,2]^4$, $f_i(a)=[0,1)$ and $f_i(*)=[0,2]$ for $i=1,2,3$. Moreover, 
$f(C^{4,*}_{4^c})=\{u^1\}$, $f(C^{4,a'}_{4^c})=\{u^2,u^3\}$, $f(C^{4,a}_{4^c})=\{u^4,u^6,u^5\}$, $f(C^{4,b'}_{4^c})=\{u^8,u^7\}$ and $f(C^{4,b}_{4^c})=\{u^9,u^{10}\}$. 
}

\medskip 

Note that Theorem \ref{kod} describes a meta-structure of a 2-periodic cube tiling of $\er^4$. To describe it we introduce some definitions. Two cubes in any cube tiling $[0,1)^d+T$ satisfy {\it Keller's condition}: For every two vectors $t,t'\in T$ there is $i\in [d]$ such that $t_i-t'_i\in \zet\setminus\{0\}$. A set $F\subset \er^d$ is said to be a {\it polycube} if there is a family of cubes $[0,1)^d+T$  satisfying Keller's condition, which is a tiling of $F$, that is if $\bigcup_{t\in T}[0,1)^d+t=F$. A non-empty polycube $B\subseteq \er^d$ is a {\it block} if it is a box  $B=B_1\times \cdots \times B_d$. 
A cube tiling   $[0,1)^d+T$ is {\it layered} if there are $a\in [0,1)$ and $i\in [d]$ such that for every $t\in T$ we have $t_i=a+z$ for some $z\in \zet$. We say that a cube tiling  $\ka T$ is {\it blockable} if there is a finite family of disjoint blocks ${\bf B}$, $|{\bf B}|>1$ such that every cube from the tiling $\ka T$ is contained in exactly one block of the family  ${\bf B}$. Clearly, every layered cube tiling of $\er^d$ is blockable. An interesting question is which non-layered cube tilings of $\er^d$ are blockable. In \cite{LP2} it was showed that every non-layered cube tiling of $\er^3$ is blockable, and in \cite{Ki} it was proved that the same is no longer true for cube tilings of $\er^d$ for $d\geq 4$. Obviously, every 2-periodic cube tiling of $\er^d$ is blockable into $2^d$ blocks, since every $f(v)$, $v\in V$, is a block (recall that $f(V)$ is a realization of a cube tiling code $V$). Therefore, regarding  a 2-periodic cube tiling of $\er^d$ it is interesting to obtain its twin pair free family of blocks ${\bf B}$ in which union of any two  blocks of ${\bf B}$ is not a block. Thus, Theorem \ref{kod} says that cubes of a 2-periodic non-layered cube tiling $\ka T$ of $\er^4$ can be arranged into one of twenty twin pair free families of blocks  ${\bf B}$. Such families can be obtained from the codes $C^i$, $i\in [20]$, given in Theorem \ref{kod} in the following way: If $\ka T$ is a  realization $f(V)$ of a cube tiling code $V$ (the way of such realization is given in Section 1), and $F_V=C^i$, then taking $A_i=f_i(a), B_i=f_i(b)$ and $f_i(*)=\er$ we have ${\bf B}=\{f_1(v_1)\times \cdots \times f_4(v_4)\colon v_1\ldots v_4\in C^i\}$. For example, if $F_V=\{a*aa, \;a**a', \;a'*a*, \;a'*a'a', \;**a'a\}$ and $f_i(a)=A$ for some fixed proper subset $A\subset \er$ and every $i\in [4]$, then 
$$
{\bf B}=\{A\times \er \times A^2, \; A\times \er^2\times A', \;  A'\times \er \times A\times \er, \; A'\times \er \times (A')^2, \; \er^2\times A' \times A \},
$$
where $A'=\er\setminus A$.

\section{Enumeration of cube tiling codes in dimension four}

Having computed all twin pair free codes $C^1-C^{20}$ together with the code $C^{21}=\{l***,l'***\}$ we are able to find the family of all codes $V\subset S^4$, where $S=\{a_1,a_1',\ldots ,a_8,a_8'\}$, that is, the family of all possible cube tiling codes in dimension four. 
To do this it is enough to replace every improper word in a code $C^i$, $i\in [21]$, by an equivalent code which contains only proper words. More precisely, if $v\in C^i$ has, up to isomorphism, the form $v=*v_2v_3v_4$, then it is replaced by two words: $u=pv_2v_3v_4$ and $w=p'v_2v_3v_4$, where $p\in S$; if, up to isomorphism, $v=**v_3v_4$, then it is replaced by four words: $u=psv_3v_4$, $w=p'sv_3v_4$,  $p=ls'v_3v_4$ and $q=l's'v_3v_4$. Since $\{***\}\equiv \{s**,s'**\}$  and  $\{***\}\equiv \{sss, s's's', *ss', s'*s, ss'*\}$,  a word $l***$ is replaced, up to isomorphism, by two words $ls**$ and $ls'**$ or by five words $lsss, ls's's', l*ss', ls'*s, lss'*$. 

After those substitutions we obtain a cube tiling code $V=V(C^i)$. About such $V$ we shall say that it is made {\it on the plane} of the code $C^i$. (Thus, Theorem \ref{kod} says that every non-layered cube tiling code is made on the plane of one of the codes $C^1-C^{20}$.) For example we find the general form of the code $V(C^{10})$. We have:
$$
\{b'**a'\}\equiv \{b'usa',b'u'sa',b'ls'a',b'l's'a'\},\;\; \{a*a'a\}\equiv \{asa'a,as'a'a\},\;\; \{a'a*a\}\equiv \{a'aqa,a'aq'a\},  
$$
$$
\{ba*a'\}\equiv \{bapa',bap'a'\},\;\; \{ba'a*\}\equiv \{ba'aq,ba'aq'\}.
$$
Thus, 
$$
V(C^{10})=\{aaaa, \; asa'a, \;as'a'a, \;a'aqa, \;a'aq'a, \;a'a'a'a, \;bapa', \;bap'a', 
$$
$$
\hspace{+22mm} ba'aq, \;ba'aq', \;ba'a'a', \;b'a'aa, \;b'lsa', \;b'l'sa', \;b'ls'a', \;b'l's'a'\},
$$
where $l,p,q,s \in S$.

\subsection{Non-isomorphic forms}

Let $S=\{a_1,a_1',\ldots ,a_8,a_8'\}$. Recall that, the elements of the group $G(S^4)$ are mappings of the form $h\circ \bar{\sigma}$, where $h$ and $\bar{\sigma}$ were defined in Section 1. It is easy to check that $G(S^d)$ is of order $4!(8!2^8)^4$. Clearly, to decide whether two codes $V,U\subset S^4$ are isomorphic it is too long to check all elements of $G(S^4)$. In this section we show which elements of $G(S^4)$ have to be considered and next we describe an algorithm that allows us to find all non-isomorphic cube tiling codes in dimension four.



Let 
$M_V=[a_{ij}]_{4\times 8}$, where $a_{ij}=(|V^{i,a_j}|,|V^{i,a'_j}|)$. We shall say that a code $V$ has a {\it compressed} form, if $M_V$ has the following two properties: For every $i\in [4]$ if $a_{ij}=(0,0)$, then $a_{ik}=(0,0)$ for every $k\geq j$. Secondly,  $|\{k\colon a_{ik}\neq (0,0)\}|\leq |\{k\colon a_{jk}\neq (0,0)\}|$ for every $i\leq j$. 
Clearly, applying step by step position bijections (compare Section 1) and next permuting positions (if needed) we may pass from $V$ into an isomorphic form $\bar{V}$ of $V$ which has a compressed form.
For example, if
$$
V=\{a_1a_3a_1a_1, a_1a_3a'_1a_1, a_1a_4a'_1a_1, a_1a'_3a_4a_1, a_1a'_3a'_4a_1, a_1a'_4a_1a'_1, a'_1a_1a_1a'_1, a'_1a_1a_3a_1,
$$
\vspace{-8mm}
$$  
\;\;\;\; \;\;\;\;\;\;a'_1a_1a'_3a_1, a'_1a'_1a_1a_2, a'_1a'_1a_1a'_2, a'_1a'_1a'_1a_1, a_2a_5a'_1a'_1, a_2a'_5a'_1a'_1, a'_2a_4a'_1a'_1,a'_2a'_4a'_1a'_1 \},
$$
then
$$
\bar{V}=\{a_1a_1a_1a_2, a_1a_1a'_1a_2, a_1a_1a'_1a_3, a_1a_1a_2a'_2, a_1a_1a'_2a'_2, a_1a'_1a_1a'_3, a'_1a'_1a_1a_1, a'_1a_1a_2a_1,  
$$
\vspace{-8mm}
$$
\;\;\;\; \;\;\;\;\;\;\;\;a'_1a_1a'_2a_1, a'_1a_2a_1a'_1, a'_1a'_2a_1a'_1, a'_1a_1a'_1a'_1, a_2a'_1a'_1a_4, a_2a'_1a'_1a'_4, a'_2a'_1a'_1a_3,a'_2a'_1a'_1a'_3 \},
$$
and

\[M_{\bar{V}}=
\begin{bmatrix}

               (6,6) & (2,2)  & (0,0)  & (0,0)  & (0,0)  & (0,0)  & (0,0)  & (0,0) \\
               (7,7) & (1,1)  & (0,0)  & (0,0)  & (0,0)  & (0,0)  & (0,0)  & (0,0) \\
               (6,6) & (1,1)  & (1,1)  & (0,0)  & (0,0)  & (0,0)  & (0,0)  & (0,0) \\
               (3,3) & (2,2)  & (2,2)  & (1,1)  & (0,0)  & (0,0)  & (0,0)  & (0,0) \\
             
\end{bmatrix}
\]

From now on, we assume that all codes have compressed forms. A permutation of rows in $M_V$, a permutation of entries in a row of $M_V$, as well as the permutation of elements in a pair $(n,n)$, where $n>0$ will be called  {\it elementary operations} on $M_V$. It is easy to see that if two codes $V,U$ are isomorphic, then the matrices $M_V$ and $M_U$ are equal, up to permutations of rows and permutations of entries in a row. (Clearly, the relation $V\approx U$ if and only if the matrices $M_V$ and $M_U$ are equal, up to permutations of rows and permutations of entries in a row is a equivalence relation.) Moreover, an isomorphism between $V$ and $U$, if it exists, may be derived from elementary operations on $M_V$. 
For example, if
  
\[M_{U}=
\begin{bmatrix}

               (1,1) & (7,7)  & (0,0)  & (0,0)  & (0,0)  & (0,0)  & (0,0)  & (0,0) \\
               (6,6) & (2,2)  & (0,0)  & (0,0)  & (0,0)  & (0,0)  & (0,0)  & (0,0) \\
               (6,6) & (1,1)  & (1,1)  & (0,0)  & (0,0)  & (0,0)  & (0,0)  & (0,0) \\
               (3,3) & (2,2)  & (2,2)  & (1,1)  & (0,0)  & (0,0)  & (0,0)  & (0,0) \\
             
\end{bmatrix}
\]
then in order to obtain $M_{\bar{V}}$ from $M_U$  we have to permute first two rows in $M_U$ and next permute the entries $(1,1)$ and $(7,7)$. Thus, applying the permutation $\sigma=(2134)$ in every word of $U$ and next applying the position bijection $(a_2,a'_2)\rightarrow (a_1,a'_1)$, $(a_1,a'_1)\rightarrow (a_2,a'_2)$ at position 2 ($(a_k,a'_k)\rightarrow (a_k,a'_k)$ for $k\in \{3,...,8\}$)  we obtain a code $U'$ which is isomorphic to $U$ and $M_{U'}=M_{\bar{V}}$.  Clearly, $U'$ need not to be isomorphic to $\bar{V}$. To check whether $U'$ and $\bar{V}$ are isomorphic we have to apply the rest mappings $h\circ \bar{\sigma}$ to $U'$, which steam from the elementary operations on $M_{U'}$ and do not change the matrix $M_{U'}$. These are: The permutation of the entries $(1,1)$ and $(1,1)$ in the third row (these induces the position bijection $(a_2,a'_2)\rightarrow (a_3,a'_3)$, $(a_3,a'_3)\rightarrow (a_2,a'_2)$ at position 3 and similarly the permutation of the entries $(2,2)$ and $(2,2)$ in the fourth row (these induces the same position bijection as above, but now at position 4). Additionally, we have to apply permutations of the form $(a_i,a'_i)\rightarrow (a'_i,a_i)$ in  every case when $|V^{j,a_i}|>0$, $j\in [4], i\in [8]$ (they steam from the permutation of elements in a pair $(n,n)$, $n>0$). Thus, to check whether $U'$ and $\bar{V}$ are isomorphic we have to check $4\times 4\times 16\times 32$ mappings $h\circ \bar{\sigma}\in G(S^4)$. 
The last factor is always the biggest as codes are compressed and thanks to the cylindrical structure of the codes we may  skip the computations that correspond to it. To do that, let $I$ consists of all  $4\times 4\times 16$ mappings $h\circ \bar{\sigma}\in G(S^4)$, $A_{\bar{V}}=\{\bar{V}_{4^c}^{4,s}\colon s\in P_{\bar{V}}\}$ and  $B_{U'}=\{(U')_{4^c}^{4,s}\colon s\in P_{U'}\}$, where $P_{U'},P_{\bar{V}}\subset S$ contain only the letters that appear at the fourth position in $U'$ and $\bar{V}$, respectively. It is easy to see that by the cylindrical structure of cube tiling code, $U'$ and $\bar{V}$ are isomorphic if and only if there is $h\circ \bar{\sigma}\in I$ such that $h\circ \bar{\sigma}(B_{U'})=A_{\bar{V}}$, where $h\circ \bar{\sigma}(B_{U'})=\{h\circ \bar{\sigma}(U')_{4^c}^{4,s}\colon s\in P_{U'}\}$.

Let $el(V,U)=1$ if $U$ can be obtained from $V$ by applying one of isomorphisms which steam from the above described elementary operations, and let $el(V,U)=0$ otherwise. Thus, $V$ and $U$ are isomorphic if and only if $el(V,U)=1$.

Before evaluating  $el(V,U)$ it is useful to check an isomorphism invariant, especially if such verification is quick. Let $V$ be a code, and let $t(V)=(t_1,t_2,t_3,t_4)$ be a vector such that $t_i$, $i\in [4]$, is a number of all twin pairs $v,u\in V$ such that $v_i=u'_i$. It is easy to see that $t(V)$ is, up to permutation, an isomorphism invariant. Let $tp(V,U)=1$ if $t(V)$ and $t(U)$ are, up to permutation, equal and $tp(V,U)=0$ otherwise. Clearly, if $tp(V,U)=0$, then $V$ and $U$ cannot be isomorphic and thus $el(V,U)=0$. 

\medskip
Let  $V\subset S^d$ be a cube tiling code (in compressed form), and let $S_i(V)$ be the set of the letters which appear in the words $v\in V$ at a $i$th position, $i\in [d]$. For example, for a code $U$ with the matrix $M_U$ given above we have  $S_1(U)=S_2(U)=\{a_1,a_1',a_2,a_2'\},S_3(U)=\{a_1,a_1',a_2,a_2',a_3,a_3'\}$ and $S_4(U)=\{a_1,a_1',...,a_4,a_4'\}$. Let $\ka N^k_4$ be a family of all non-isomorphic cube tiling codes in dimension four such that $S_4(V)=\{a_1,a_1',...,a_k,a_k'\}$. To compute a set $\ka N_4$ of all non-isomorphic cube tiling codes in dimension four we shall divide the computations into three parts: In the first part we shall indicate a family $\ka N_4^6$. To do this, we shall proceed as follows: 


\medskip
{\bf Algorithm 2}.

\medskip
{\it Input}. The family $\{C^i\colon i\in [21]\}$.

{\it Output}. The family $\ka N^6_4$.

\smallskip
1. Compute all codes $V\subset S^4$, $S=\{a_1,a'_1,...,a_6,a'_6\}$, on the plane of the codes $C^i$ for $i \in [21]$ in the manner described in the previous section. Next, every such code is transformed into compressed form. Denote the family of all such codes by $\ka F$. 

\smallskip
2. Divide the set $\ka F$ into equivalence classes $(\ka F_i)_{i\in [m]}$ of the relation $\approx$. ($V\approx U$ if and only if the matrices $M_V$ and $M_U$ are equal, up to permutations of rows and permutations of entries in a row.)


\smallskip
3. For every $i\in [m]$ indicate the set $\ka N^6_4(i)$ consisting of all non-isomorphic codes in $\ka F_i$. To find out whether codes $V,U \in \ka F_i$ are isomorphic  evaluate first $tp(V,U)$ and if $tp(V,U)=1$, evaluate $el(V,U)$.  
Since codes from different sets $\ka N^6_4(i)$ and $\ka N^6_4(j)$ cannot be isomorphic, we have $\ka N^6_4=\bigcup_{i \in [m]} \ka N^6_4(i)$.

\bigskip
Since indication of $\ka N_4^7$ and $\ka N_4^8$ according the above algorithm would be unnecessarily long we shall take a shortcut to compute these two families.
Let $\ka N_3$ be a family of all non-isomorphic cube tilings in dimension three. (The family $\ka N_3$ can be easily indicated even by hand computations.) We have $|\ka N_3|=17$. The set $\ka N_4^8$ can be derived immediately from $\ka N_3$ thanks to the cylindrical structure of a cube tiling code. To see this, let $V\subset S^{d-1}$ be a cube tiling code and let $\pi^k_V=\{V^1,...,V^k\}$ be a partition of $V$. Let  $Vs=\{v_1...v_{d-1}s\colon v\in V\}$ for $s\in S$ . It follows from the cylindrical structure of a cube tiling code that for every pairwise different letters $\{a_{i_1},a_{i_1}'...,a_{i_k},a_{i_k}'\}\subset S$ and every cube tiling codes $V,W\subset S^{d-1}$ the set of words
\begin{equation}
\label{ws}
U(\pi^k_V,\pi^k_W)=V^1a_{i_1}\cup W^1a_{i_1}'\cup\ldots \cup V^ka_{i_k}\cup W^ka_{i_k}'  
\end{equation}
is a cube tiling code in dimension $d$ if and only if $V^i\equiv W^i$ for every $i\in [k]$ (compare \cite{GKP,Ost}). Since $S_4(V)=\{a_1,a_1',...,a_8,a_8'\}$ for $V\in \ka N_4^8$, we have $V^i=W^i$ for every $i\in [8]$. Thus, 
$$
\ka N_4^8=\{U(\pi^8_V,\pi^8_V)\colon V\in \ka N_3\}.
$$
In the case of the set $\ka N_4^7$ we consider a family of  codes $U(\pi^7_V,\pi^7_W)$ where $V,W$ are cube tiling codes in dimension three. Since $|V|=8$,  we have $|V^i|=|W^i|=2$ for exactly one $i\in [7]$. It is easy to check that in this case $V^i\equiv W^i$ and  $V^i\neq W^i$ if and only if both codes $V^i, W^i$ are twin pairs. Let $\ka A$ be the family of the codes $U(\pi^7_V,\pi^7_V)$, $V\in \ka N_3$, such that $V^i$ is not a twin pair for each $i\in [7]$. Moreover, let $\ka B$ the family of all codes $U(\pi^7_V,\pi^7_W)$, $V\in \ka N_3$, where $V^i$ is a twin pair for some $i\in [7]$ and $W$ is a cube tiling code in dimension three such that $W^i\equiv V^i$ and $W^j=V^j$ for $j\neq i$.  Since two codes $V\in \ka A$ and $U\in \ka B$ cannot be isomorphic, to obtain $\ka N_4^7$ it is enough to compute sets $\ka N(\ka A)$ and $\ka N(\ka B)$ of all non-isomorphic codes in $\ka A$ and $\ka B$, respectively. We do it in the manner described in the second and the third step  of Algorithm 2, where $\ka F\in \{\ka A,\ka B\}$. After this, we get
$$
\ka N^7_4=\ka N(\ka A)\cup \ka N(\ka B)
$$
and finally
$$
\ka N_4=\ka N^6_4 \cup \ka N^7_4\cup \ka N^8_4.
$$

\subsection{The number of all cube tiling codes in dimension four}
To compute the number of all cube tiling codes $V\subset S^d$ we may first find some smaller set of codes. To do this, let us note that  for a given code $V\subset S^d$  the previously defined set of the letters $S_i(V)$ is usually a proper subset of $S$. 
For $i \in [d]$ let $H_i(V)$ be the set of position bijections $h_i$ such that $h_i(S_i(V))=S_i(V)$. Moreover, let $\Sigma(V)$ consists of all permutations $\sigma$ of $[d]$ such that if $\sigma(i)\neq i$, then $S_i(V)=S_{\sigma(i)}(V)$ for $i \in [d]$. Observe that if $h\circ \bar{\sigma}=(h_1,...,h_d)\circ \bar{\sigma}\in G(S^d)$ is such that $h_i\not \in H_i(V)$ for some $i\in [d]$ or $\sigma \not \in \Sigma(V)$ (recall that $\bar{\sigma}(v)=v_{\sigma(1)}\ldots v_{\sigma(d)}$), then $V\neq h\circ \bar{\sigma}(V)$. Let us call the set 
$$
o_m(V)=\{h\circ \bar{\sigma}(V)\colon \sigma \in \Sigma(V) \; {\rm and}\; h_i\in H_i\; {\rm for}\; i\in [d]\}
$$
a {\it minimal orbit} of $V$.  
It is easy to see that that having the cardinality of the minimal orbit of a cube tiling code $V\subset S^d$ we are able to calculate the cardinality of the whole orbit $o(V)$, since
\begin{equation}
\label{wz}
|o(V)|= (d!/|\Sigma(V)|){k\choose k_1}\cdots {k\choose k_d} |o_{m}(V)|,
\end{equation}
where $S=\{a_1,a_1',...,a_k,a_k'\}$ and $k_i=\frac{1}{2}|S_i(V)|$ for $i\in [d]$. Therefore, to find the number of all cube tiling codes in dimension four, which is equal to $M_d=\sum_{V\in \ka N_d}|o(V)|$, where $\ka N_d$ is a family of all non-isomorphic cube tiling codes in dimension $d$, it is enough to compute the numbers $(|o_m(V)|)_{V\in \ka N_d}$ 
and next apply the formula (\ref{wz}). It seems that beside $M_d$ it is worth to find the sum  $M^m_d=\sum_{V\in \ka N_d}|o_m(V)|$. We call $M^m_d$ the {\it total number of minimal orbits}.


\begin{uw}
\label{orb}
{\rm
If the number $|o_m(V)|$ is interpreted as a measure of the complexity of the code $V$ it should be low for the simplest codes. A cube tiling code $V$ whose combinatorics is the simplest 
is called {\it simple} (or {\it regular}) code in which $v_i=u_i$ or $v_i=u_i'$ for every $v,u \in V$ and every $i\in [d]$ (clearly, it is isomorphic to the binary code $\{0,1\}^d$). For example, $V=\{aa,aa',a'a,a'a'\}$ is a simple cube tiling code in dimension two. Then $S_i(V)=\{a,a'\}$ for $i=1,2$. If we consider the family of all cube tiling codes in dimension two we have to, by Lemma \ref{uu}, consider codes written down in the alphabet $S=\{a,a',b,b'\}$. Then $|o(V)|=4$ but $|o_{m}(V)|=1$. Similarly, for a simple code $V\subset S^4$, $S=\{a_1,a'_1,...,a_8,a'_8\}$, we have $|o(V)|=4096$ and $|o_{m}(V)|=1$.

On the other hand every cube tiling code is build up from pieces of simple codes called {\it simple components}, that is, $V=P^1\cup\ldots \cup P^n$, where $P^i=V\cap U^i\neq\emptyset$ and $U^i$ is a simple cube tiling code for $i\in [n]$. (The structure of simple components $P^i$ are subjected certain interesting restrictions, see \cite{KP}.) In this approach, for the simple code $V\in \ka N_d$ the number $|o(V)|$ says how many `building blocks' to create cube tiling codes we have. Therefore, in the next section we give both  numbers $M_4$ and $M_4^{m}$.

}
\end{uw}

\subsection{Results}
In the computations we used 8-core 3.4-GHz processor and the indication of a set $\ka N_4$ consisting of all non-isomorphic cube tiling codes in dimension four took about two days (The code was written in Python.) By Lemma \ref{uu}, to obtain these cube tiling codes it is enough to consider codes that are written down in the alphabet $S=\{a_1,a'_1,...,a_8,a'_8\}$.

\begin{tw}
\label{num}
There are $N_4=27,385$ non-isomorphic cube tiling codes in dimension four, and the total number of such codes is equal to 
$M_4=17,794,836,080,455,680$. Moreover, $M^{m}_4=1,108,646,656$. In particular, there are $N_4$ non-isomorphic $2$-periodic cube tilings of $\er^4$, cube tilings of $\te^4$ and $r$-perfect codes in $\zet_{4r+2}$ in the maximum metric, and the total number and the total number of minimal orbits of such tilings and codes are equal to $M_4$ and $M_4^m$, respectively. 

\hfill{$\square$}
\end{tw}


Below we present the cardinalities $o$ of orbits $o(V)$ of the codes $V\in \ka N_4$ and  the number $N(o)$ of codes  in $\ka N_4$ with the cardinality $o$. 


\begin{center}
\begin{tabular}{ c c|c c| c c| c c }
\hline
$o$   &  $N(o)$ & $o$ & $N(o)$ & $o$   &  $N(o)$ & $o$ & $N(o)$\\

\hline

4,096 & 1 & 924,844,032 & 105 & 37,764,464,640 & 2 & 906,347,151,360 & 45\\
 114,688 & 1 & 944,111,616 & 4 & 38,843,449,344 & 1091 & 932,242,784,256 & 52\\
 344,064 & 2 & 990,904,320 & 9 & 41,617,981,440 & 23 & 971,086,233,600 & 426\\
 458,752 & 2 & 1,011,548,160 & 2 & 45,317,357,568 & 270 & 1,019,640,545,280 & 5\\
 917,504 & 2 & 1,078,984,704 & 128 & 48,554,311,680 & 282 & 1,109,812,838,400 & 26\\
 1,204,224 & 1 & 1,083,801,600 & 1 & 52,022,476,800 & 2 & 1,165,303,480,320 & 493\\
 1,376,256 & 4 & 1,156,055,040 & 286 & 55,490,641,920 & 293 & 1,223,568,654,336 & 4\\
 2,408,448 & 2 & 1,213,857,792 & 4 & 56,646,696,960 & 29 & 1,359,520,727,040 & 136\\
 2,752,512 & 4 & 1,258,815,488 & 1 & 58,265,174,016 & 83 & 1,456,629,350,400 & 42\\
 3,440,640 & 1 & 1,387,266,048 & 184 & 60,692,889,600 & 1 & 1,553,737,973,760 & 298\\
 4,128,768 & 5 & 1,618,477,056 & 235 & 64,739,082,240 & 594 & 1,631,424,872,448 & 29\\
 4,816,896 & 5 & 1,734,082,560 & 27 & 67,976,036,352 & 52 & 1,664,719,257,600 & 11\\
 5,505,024 & 2 & 1,849,688,064 & 5 & 69,363,302,400 & 213 & 1,699,400,908,800 & 2\\
 5,619,712 & 1 & 1,888,223,232 & 35 & 72,831,467,520 & 1 & 1,747,955,220,480 & 5\\
 8,257,536 & 6 & 1,981,808,640 & 1 & 75,528,929,280 & 3 & 1,812,694,302,720 & 18\\
 9,633,792 & 11 & 2,023,096,320 & 7 & 77,686,898,688 & 738 & 1,942,172,467,200 & 310\\
 10,321,920 & 2 & 2,157,969,408 & 27 & 83,235,962,880 & 56 & 2,039,281,090,560 & 15\\
 13,762,560 & 1 & 2,312,110,080 & 476 & 90,634,715,136 & 89 & 2,219,625,676,800 & 5\\
 16,515,072 & 23 & 2,427,715,584 & 10 & 97,108,623,360 & 1052 & 2,330,606,960,640 & 399\\
 19,267,584 & 25 & 2,517,630,976 & 2 & 101,964,054,528 & 1 & 2,447,137,308,672 & 3\\
 20,643,840 & 5 & 2,774,532,096 & 267 & 104,044,953,600 & 4 & 2,719,041,454,080 & 70\\
 27,525,120 & 3 & 3,236,954,112 & 758 & 110,981,283,840 & 168 & 2,913,258,700,800 & 136\\
 28,901,376 & 3 & 3,468,165,120 & 112 & 113,293,393,920 & 96 & 3,107,475,947,520 & 18\\
 33,030,144 & 6 & 3,776,446,464 & 58 & 116,530,348,032 & 171 & 3,262,849,744,896 & 6\\
 33,718,272 & 3 & 4,046,192,640 & 34 & 121,385,779,200 & 16 & 3,329,438,515,200 & 1\\
 38,535,168 & 50 & 4,335,206,400 & 1 & 129,478,164,480 & 416 & 3,398,801,817,600 & 7\\
 41,287,680 & 9 & 4,624,220,160 & 283 & 135,952,072,704 & 147 & 3,495,910,440,960 & 10\\
 44,957,696 & 2 & 4,855,431,168 & 47 & 138,726,604,800 & 182 & 3,641,573,376,000 & 2\\
 55,050,240 & 1 & 5,549,064,192 & 249 & 145,662,935,040 & 28 & 3,884,344,934,400 & 218\\
 57,802,752 & 27 & 5,664,669,696 & 10 & 155,373,797,376 & 88 & 4,078,562,181,120 & 32\\
 67,436,544 & 5 & 6,473,908,224 & 962 & 166,471,925,760 & 53 & 4,369,888,051,200 & 1\\
 72,253,440 & 3 & 6,936,330,240 & 351 & 169,940,090,880 & 2 & 4,661,213,921,280 & 174\\
 77,070,336 & 28 & 7,283,146,752 & 1 & 194,217,246,720 & 1557 & 5,438,082,908,160 & 20\\
 78,675,968 & 1 & 7,552,892,928 & 69 & 203,928,109,056 & 9 & 5,826,517,401,600 & 175\\
 82,575,360 & 31 & 8,092,385,280 & 188 & 208,089,907,200 & 17 & 6,117,843,271,680 & 1

\end{tabular}
\end{center} 
 
\begin{center}
\begin{tabular}{ c c|c c| c c| c c }
\hline
$o$   &  $N(o)$ & $o$ & $N(o)$ & $o$   &  $N(o)$ & $o$ & $N(o)$\\

\hline
 86,704,128 & 4 & 8,670,412,800 & 13 & 221,962,567,680 & 13 & 6,797,603,635,200 & 5\\
 89,915,392 & 2 & 9,248,440,320 & 223 & 226,586,787,840 & 152 & 6,991,820,881,920 & 19\\
 110,100,480 & 3 & 9,441,116,160 & 1 & 233,060,696,064 & 259 & 7,283,146,752,000 & 6\\
 115,605,504 & 99 & 9,710,862,336 & 252 & 242,771,558,400 & 76 & 7,768,689,868,800 & 99\\
 134,873,088 & 19 & 10,404,495,360 & 4 & 258,956,328,960 & 70 & 8,157,124,362,240 & 24\\
 144,506,880 & 6 & 11,098,128,384 & 34 & 271,904,145,408 & 186 & 8,739,776,102,400 & 2\\
 154,140,672 & 8 & 11,329,339,392 & 63 & 277,453,209,600 & 155 & 9,322,427,842,560 & 54\\
 165,150,720 & 26 & 12,138,577,920 & 5 & 291,325,870,080 & 116 & 9,710,862,336,000 & 1\\
 173,408,256 & 8 & 12,947,816,448 & 311 & 302,115,717,120 & 3 & 10,196,405,452,800 & 1\\
 202,309,632 & 1 & 13,872,660,480 & 577 & 312,134,860,800 & 1 & 10,876,165,816,320 & 1\\
 231,211,008 & 227 & 14,566,293,504 & 3 & 332,943,851,520 & 18 & 11,653,034,803,200 & 126\\
 247,726,080 & 5 & 15,105,785,856 & 12 & 339,880,181,760 & 24 & 12,235,686,543,360 & 2\\
 269,746,176 & 74 & 15,173,222,400 & 1 & 364,157,337,600 & 1 & 13,595,207,270,400 & 1\\
 289,013,760 & 21 & 16,184,770,560 & 703 & 388,434,493,440 & 1372 & 13,983,641,763,840 & 12\\
 314,703,872 & 1 & 16,994,009,088 & 2 & 407,856,218,112 & 24 & 14,566,293,504,000 & 17\\
 330,301,440 & 20 & 17,340,825,600 & 36 & 416,179,814,400 & 26 & 15,537,379,737,600 & 26\\
 346,816,512 & 20 & 18,496,880,640 & 37 & 453,173,575,680 & 101 & 16,314,248,724,480 & 3\\
 359,661,568 & 1 & 19,421,724,672 & 747 & 466,121,392,128 & 219 & 17,479,552,204,800 & 7\\
 404,619,264 & 5 & 20,808,990,720 & 8 & 485,543,116,800 & 266 & 18,644,855,685,120 & 7\\
 462,422,016 & 302 & 22,658,678,784 & 176 & 543,808,290,816 & 70 & 20,392,810,905,600 & 2\\
 472,055,808 & 1 & 24,277,155,840 & 34 & 554,906,419,200 & 86 & 23,306,069,606,400 & 55\\
 495,452,160 & 1 & 25,895,632,896 & 17 & 582,651,740,160 & 297 & 27,967,283,527,680 & 3\\
 505,774,080 & 1 & 27,745,320,960 & 605 & 604,231,434,240 & 1 & 29,132,587,008,000 & 14\\
 539,492,352 & 155 & 28,323,348,480 & 4 & 624,269,721,600 & 1 & 31,074,759,475,200 & 5\\
 578,027,520 & 106 & 29,132,587,008 & 26 & 679,760,363,520 & 115 & 34,959,104,409,600 & 5\\
 629,407,744 & 3 & 30,211,571,712 & 1 & 728,314,675,200 & 7 & 43,698,880,512,000 & 1\\
 660,602,880 & 10 & 32,369,541,120 & 1060 & 776,868,986,880 & 580 & 46,612,139,212,800 & 14\\
 693,633,024 & 63 & 33,988,018,176 & 8 & 815,712,436,224 & 38 & 58,265,174,016,000 & 8\\
 809,238,528 & 46 & 34,681,651,200 & 117 & 832,359,628,800 & 22 & 69,918,208,819,200 & 1\\ 
 867,041,280 & 3 & 36,415,733,760 & 1 & 873,977,610,240 & 2  
  
\end{tabular}
\end{center}
\noindent{\footnotesize Table 1: The cardinalities $o$ and  the number $N(o)$ of codes  in $\ka N_4$ with the cardinality $o$.
}

\medskip

\medskip
In Table 2 we present the numbers $n(c)$ of codes $V \in \ka N_4$ with the distribution $c=(c_1,c_2,c_3,c_4)$, where $c_i$, $i\in [4]$, is the number of different pairs of letters $(a_j,a'_j)$, $j\in [8]$, which appear in the code $V$ at the $i$th position, that is, $c_i=|\{j\in [8]\colon V^{i,a_j}\cup V^{i,a'_j}\neq \emptyset\}|$. In more geometric interpretation, $c_i$ is the number of cylinders in the $i$th direction in a realization $f(V)$ (compare Figure 3). 

\pagebreak
We give this classification as the vector $c$ reflects the structure of a cube tiling code. 
The codes are presented in compressed forms. There are 90 different distributions $c$.

\begin{center}
\begin{tabular}{ c c| c c| c c}
\hline
   $c$   &  $n(c)$  &   $c$   &  $n(c)$ &   $c$   &  $n(c)$\\

\hline

$(1, 1, 1, 1)_b$ & 1 &  (1, 2, 2, 2) & 417 &  (1, 3, 4, 6) & 12\\
 (1, 1, 1, 2) & 13 &  (1, 2, 2, 3) & 2335 &  $(1, 3, 4, 7)_l$ & 1\\
 (1, 1, 1, 3) & 44 &  (1, 2, 2, 4) & 2120 &  (1, 3, 5, 5) & 9\\
 (1, 1, 1, 4) & 66 &  (1, 2, 2, 5) & 955 &  $(1, 3, 5, 6)_l$ & 1\\
 (1, 1, 1, 5) & 43 &  (1, 2, 2, 6) & 236 &  (1, 4, 4, 4) & 32\\
 (1, 1, 1, 6) & 15 &  (1, 2, 2, 7) & 32 &  (1, 4, 4, 5) & 14\\
 (1, 1, 1, 7) & 3 &  (1, 2, 2, 8) & 4 &  $(1, 4, 4, 6)_l$ & 1\\
 (1, 1, 1, 8) & 1 &  (1, 2, 3, 3) & 2972 &  $(1, 4, 5, 5)_l$ & 2\\
 (1, 1, 2, 2) & 130 &  (1, 2, 3, 4) & 3656 &  $(2, 2, 2, 2)_b$ & 183\\
 (1, 1, 2, 3) & 768 &  (1, 2, 3, 5) & 1270 &  (2, 2, 2, 3) & 613\\
 (1, 1, 2, 4) & 938 &  (1, 2, 3, 6) & 257 &  (2, 2, 2, 4) & 332\\
 (1, 1, 2, 5) & 503 &  (1, 2, 3, 7) & 33 &  (2, 2, 2, 5) & 120\\
 (1, 1, 2, 6) & 141 &  (1, 2, 3, 8) & 3 &  (2, 2, 2, 6) & 29\\
 (1, 1, 2, 7) & 21 &  (1, 2, 4, 4) & 906 &  (2, 2, 2, 7) & 5\\
 (1, 1, 2, 8) & 3 &  (1, 2, 4, 5) & 478 &  (2, 2, 2, 8) & 1\\
 (1, 1, 3, 3) & 829 &  (1, 2, 4, 6) & 79 &  (2, 2, 3, 3) & 461\\
 (1, 1, 3, 4) & 1393 &  (1, 2, 4, 7) & 8 &  (2, 2, 3, 4) & 273\\
 (1, 1, 3, 5) & 591 &  $(1, 2, 4, 8)_l$ & 1 &  (2, 2, 3, 5) & 63\\
 (1, 1, 3, 6) & 140 &  (1, 2, 5, 5) & 60 &  (2, 2, 3, 6) & 9\\
 (1, 1, 3, 7) & 20 &  (1, 2, 5, 6) & 14 &  (2, 2, 3, 7) & 1\\
 (1, 1, 3, 8) & 2 &  $(1, 2, 5, 7)_l$ & 1 &  (2, 2, 4, 4) & 14\\
 (1, 1, 4, 4) & 479 &  $(1, 2, 6, 6)_l$ & 2 &  (2, 2, 4, 5) & 1\\
 (1, 1, 4, 5) & 309 &  (1, 3, 3, 3) & 819 &  (2, 3, 3, 3) & 109\\
 (1, 1, 4, 6) & 61 &  (1, 3, 3, 4) & 1013 &  (2, 3, 3, 4) & 57\\
 (1, 1, 4, 7) & 7 &  (1, 3, 3, 5) & 264 &  (2, 3, 3, 5) & 10\\
 (1, 1, 4, 8) & 1 &  (1, 3, 3, 6) & 44 &  (2, 3, 3, 6) & 1\\
 (1, 1, 5, 5) & 47 &  (1, 3, 3, 7) & 6 &  (2, 3, 4, 4) & 3\\
 (1, 1, 5, 6) & 13 &  $(1, 3, 3, 8)_l$ & 1 &  $(3, 3, 3, 3)_b$ & 9\\
 (1, 1, 5, 7) & 1 &  (1, 3, 4, 4) & 325 &  (3, 3, 3, 4) & 3\\
 (1, 1, 6, 6) & 2 &  (1, 3, 4, 5) & 119 &  (3, 3, 3, 5) & 1

 \end{tabular}
\end{center}

\noindent{\footnotesize Table 2: The distributions of cylinders. A lamination is marked by $c_l$, and a balanced code by $c_b$.
}

\medskip

There are 10 {\it laminations}, the codes with the greatest possible number of cylinders in all directions. It was conjecture in \cite{DI} that this number is equal to $2^d-1$ and confirmed in \cite{Pr}.  Moreover, there are 193 {\it balanced} cube tiling codes (named in \cite{Ba} as {\it $EDL$-partition}), that is, codes in which the number of cylinders in all directions is the same. 

In the next table we present a distribution of codes $V\in \ka N_4$ depending on a number of letters used in $V$. For $k\in [8]$ let $\ka N_4^k$ consists of all codes $V \in \ka N_4$ such that $V\subset \{a_1,a_1',...,a_k,a_k'\}^4$ and there is $ i \in [4]$ with $V^{i,a_j}\cup V^{i,a_j'}\neq \emptyset$ for every $j\in [k]$. In other words, $V\in \ka N_4^k$ if and only if $S_4(V)=\{a_1,a_1',...,a_k,a_k'\}$ (compare the previous section). We have

\begin{center}
\begin{tabular}{ c c| c c| c c|cc}
\hline
   $k$   &  $|\ka N_4^k|$  &   $k$   &  $|\ka N_4^k|$ &   $k$   &  $|\ka N_4^k|$ & $k$   &  $|\ka N_4^k|$\\

\hline
1 &  1 & 3 &  8959 & 5 &  4859 & 7 &  139\\
2 &  744 & 4 &  11610 & 6 &  1057 & 8 &  17
\end{tabular}
\end{center}

\noindent{\footnotesize Table 3: The cardinalities of the sets $\ka N_4^k$ for $k\in [8]$.
}

\medskip

Finally we present the current state of the classifications of cube tilings codes:

\begin{center}
\begin{tabular}{  c|  c|  c}
\hline
   $d$   &  $N^2_d$  &   $N_d$\\

\hline
$d=1$ &  1 & 1 \\
$d=2$ &  2 & 2 \\
$d=3$ &  9 & 17 \\
$d=4$ &  744 & 27,385 \\
$d=5$ &  899,710,227 & $?$\\
$d\geq 6$ & $?$ & $?$ 
\end{tabular}
\end{center}

\noindent{\footnotesize Table 4: In the second column the number $N^2_d$ of non-isomorphic cube tiling codes in dimension $d$ that are written down in the alphabet $S=\{a,a',b,b'\}$ is given. The number $N_d$ of all non-isomorphic cube tiling codes is given in the last column. 
}

\subsection{A note on the number of cube tiling codes in dimension five}
At the end of this section let us discus briefly the problem of an estimation of the number $M_5$. To do this, we shall give a general formula for  $M_d=|\ka M_d|$, where $\ka M_d$ is the set of all cube tiling codes in dimension $d$ written down in the alphabet $S=\{a_1,a_1',...,a_{2^{d-1}},a_{2^{d-1}}'\}$. 
Let $S(n,k)$ be the number of all partitions of an $n$-element set into  $k$ non-empty subsets (the Stirling number of the second kind), and let $C_n=\sum_{k=2}^{n}$${n}\choose {k}$$S(n,k)k!$.  Since the codes from $\ka M_d$ are written down in the alphabet $S$ having $2^{d-1}$ pairs of letters $a_i,a'_i$, the number of all layered cube tiling codes $U(\pi^1_V,\pi^1_W) \subset S^d$ (compare (\ref{ws})) is equal to $2^{d-1}(M_{d-1}(S))^2$, where $M_{d-1}(S)$ is the number of all cube tiling codes in dimension $d-1$ but written down in the alphabet $S$. Additionally, the number of all codes $U(\pi^k_V,\pi^k_V)\in  \ka M_d$ for $k\in \{2,3,...,2^{d-1}\}$ is equal to $M_{d-1}(S)C_{2^{d-1}}$. Finally, let $\ka U_d$ be the family of all cube tiling codes $U\subset S^d$ such that  $U=U(\pi^k_V,\pi^k_W)$, where  $k\in \{2,3,...,2^{d-1}\}$ and $V^i\neq W^i$ for at least one $i\in [k]$. Thus,

\begin{equation}
\label{li}
M_d=2^{d-1}(M_{d-1}(S))^2+M_{d-1}(S)C_{2^{d-1}}+|\ka U_d|.
\end{equation}

Unfortunately, we are not able to indicate the cardinality of $\ka U_d$ based only on the numbers $M_{d-1}$ and $C_{2^{d-1}}$. The number $|\ka U_d|$ depends on pairs of polybox codes $V_i,W_i\subset S^{d-1}$ which are equivalent and $V_i\neq W_i$. The structure of some elements in $\ka U_d$ is substantially new and cannot be realized by any cube tiling code in lower dimensions. 
For example, all cube tiling codes in dimension two are layered, but in dimension three there is an element in $\ka U_3$ which is not a layered code (Figure 2). Codes in $\ka U_d$ depend on rigid polybox codes: A polybox code $Q$ is {\it rigid} if $Q\equiv P$ implies $Q=P$. For  $U(\pi^k_V,\pi^k_W)\in \ka U_d$ at least one code  $V^i$ for some $i\in [k]$ is not rigid. It is known (\cite{KisL}) that for $d=4$ there is, up to isomorphism, only one twin pair free polybox code $Q$ which is not rigid ($Q$ contains twelve words). All other non-rigid polybox codes $Q\subset S^4$ have to contain a twin pair.  So, let us suppose that in a partition $\pi^k_V$ of a cube tiling code $V\subset S^{d-1}$ there are $m\leq k$ codes $V^i$ from the representation (\ref{ws}) each containing a twin pair. Observe now, that if $Q=\{v,u\}$ is a twin pair in which for example $v_r=a_1$ and $u_r'=a_1'$, then twin pairs $Q^j=\{v^j,u^j\}$, where $v_r^j=a_j,u_r^j=a_j'$, $j \in \{2,...,16\}$ and $v_n=u_n=v^j_n=u^j_n$ for $n\in [4]\setminus\{r\}$, are equivalent polybox codes to $Q$ and $Q\cap Q^j=\emptyset$   for all $j$. Therefore, if $Q\subset V^i$, then there are at least $15$ different codes $W^i$ such that $V^i\equiv W^i$ and $V^i\neq W^i$, as it is enough to take $W^i=V^i\setminus Q \cup Q^j$. Thus, such $V$ generates at least $(16^{m}-1)$${16}\choose {k}$$(2^k-2)m!(k-m)!$ codes of the form $U(\pi^k_V,\pi^k_W)\in \ka U_5$. It seems that using Theorem \ref{kod} it is possible to give some estimation of $|\ka U_5|$. We will not however try to do it here. Instead of it, we give the sum $\tilde{M}_5=2^{4}(M_{4}(S))^2+M_{4}(S)C_{2^{4}}$. We have $M_4(S)=619,671,688,833,358,364,672\;\;{\rm and}\;\; C_{16}=18,446,723,150,919,663,616$ and consequently $\tilde{M}_5=6,155,318,955,547,189,794,969,842,157,424,520,300,855,296.$ Thus, 
$$
M_5>6\times 10^{42}.
$$

\section{Glue and cut procedure}

Recall that a twin pair $v,u$ is glued at the $i$th position if $v,u$ is replaced by the improper word $w$ having  the star at the $i$th position where $v_i=u_i'$ and $w_j=v_j$ for all $j\neq i$ (we called the word $w$ a gluing of $v$ and $u$). Conversely, if an improper word $w$ with exactly one star at the $i$th position is replaced by a twin pair $v,u$ such that $v_i=s,u_i=s'$, where $s\neq *$ and $v_j=u_j=w_j$ for all $j\neq i$, then we say that $w$ is {\it cut} at the $i$th position. Let $V,U\subset (*S)^d$ be equivalent polybox codes. We say that the polybox code $U$ is obtained  from the polybox code $V$  by the {\it glue and cut procedure} if there is a sequence of gluing and cutting leading from $V$ to $U$. (The glue and cut procedure was considered by K. Przes\l awski and the author independently from Dutour and Itoh (\cite{DI2}). In \cite{DI2,Ost} this procedure is called {\it switching}.) 
For example, 
in Figure 5 we see a visualization of passing from one polybox code to another. 

In \cite{DI2} it was shown that for all $d\leq 4$ we may pass from an arbitrary cube tiling code $V\subset \{a,a',b,b'\}^d$ to an arbitrary cube tiling code $U\subset \{a,a',b,b'\}^d$ by the glue and cut procedure, and in \cite{Ost} it was proved that the same is true for any two cube tiling codes $V,U\subset \{a,a',b,b'\}^5$.

Using the above result we are able to prove 

\begin{tw}
\label{cut}
For every $d\leq 5$ and $S$ it is possible to pass from a cube tiling code $V\subset S^d$ to a cube tiling code $U\subset S^d$ by the glue and cut procedure.
\end{tw}
\proof
By the results given before the theorem, it is enough to show that one may pass by the glue and cut procedure form $V$ to a cube tiling code $\bar{V}\subset \{a,a',b,b'\}^d$.
It is easy to show that this is possible for $d\leq 3$.

Let $d=4$. A twin pair free code $F_V$ of $V$  has one of the forms given in Theorem \ref{kod} or $V^{n-1}$ is, up to isomorphism, of the form $\{a***,a'***\}$ (see Section 2 for $V^{n-1}$).
Thus, every word $v\in F_V\cup V^{n-1}$ contains at most three stars. If $v$ contains one, two or three stars and $U_v\subset V$ is the set of words that has been aggregated to $v$, that is, $\{v\}\equiv U_v$, then we may identify $U_v$ with  a cube tiling code in dimension one, two or three, respectively. Thus, by the result mentioned above, we may pass by the glue and cut procedure, from $U_v$ to a code $\bar{U}_v\subset \{a,a',b,b'\}^4$. Since, up to isomorphism, $F_V,V^{n-1}\subset \{a,a',b,b',*\}^4$, it follows that we may pass by the glue and cut procedure from $V$ to a cube tiling code  $\bar{V}\subset \{a,a',b,b'\}^4$.


Let $d=5$. 
If $F_V=\{*****\}$, then, up to isomorphism, $V^{n-1}=\{a****,a'****\}$. As above, for $v\in V^{n-1}$ we may identify the codes $U_v\subset V$, where $\{v\}\equiv U_v$ with cube tiling codes in dimension four.  Thus, by the first part of the proof we may pass by the glue and cut procedure from $V$ to  $\bar{V}\subset \{a,a',b,b'\}^5$. Clearly, the same is true if $F_V\neq \{*****\}$ and $F_V\subset \{a,a',b,b',*\}^5$, as each word in $F_V$ contains at most four stars.

 
So assume that there is $i\in [5]$, say $i=1$,  and there are $a,b,c\in S$ such that $F_V^{1,s}\cup F_V^{1,s'}\neq \emptyset$ for $s\in \{a,b,c\}$. Thus, for at least one $s$, say it is $s=c$, we have $\sum_{v\in F_V^{1,s}}|v|\leq 5$ for $s \in \{c,c'\}$ (see (\ref{nor})). Since $F_V$ is a twin pair free code, it is not hard to compute that the set $F_V^{1,c}\cup F_V^{1,c'}$ has, up to isomorphism, one of the following forms:  

\begin{enumerate}
\item $F_V^{1,c}=\{ca**a, caaaa'\}\;\; {\rm and}\;\; F_V^{1,c'}=\{c'a*a'a, c'aaa*,   c'aa'aa\}$
\item $F_V^{1,c}=\{ca*a'a, ca*ba',  caab'a'\}\;\; {\rm and}\;\;F_V^{1,c}=\{c'aaa'*, c'aaaa',  c'aa'a'a, c'aa'ba'\}$
\item $F_V^{1,c}=\{ ca*aa, caa*a', caa'a'a'\}\;\; {\rm and}\;\;  F_V^{1,c'}=\{c'a*a'a', c'aaa*,  c'aa'aa\}$
\item $F_V^{1,c}=\{c*aaa, ca*a'b, caaaa'\}\;\; {\rm and}\;\;  F_V^{1,c'}=\{c'aa*b, c'aaab', c'aa'a'b, c'a'aaa\}$
\item $F_V^{1,c}=\{ca*aa, caa*a'\}\;\; {\rm and}\;\;  F_V^{1,c'}=\{c'aaa*, c'aa'aa, c'aaa'a'\}$
\item $F_V^{1,c}=\{ca*aa, caaa'a\}\;\; {\rm and}\;\;  F_V^{1,c'}=\{c'aa*a, c'aa'aa\}$
\end{enumerate}
Let $U\subset V^{1,c}\cup V^{1,c'}$ be such that $U\equiv F_V^{1,c}\cup F_V^{1,c'}$ (that is, $U$ consists of words $v\in V$ which have been glued to words from $F_V^{1,c}\cup F_V^{1,c'}$). Knowing the forms of $F_V^{1,c}\cup F_V^{1,c'}$ we may check that it is possible to pass from $U$ to a polybox code $W$ by the glue and cut procedure, where $W$ contains only words with $b$ or $b'$ at the first position. To do that, note that $U^{1,c}\equiv F_V^{1,c}$ and  $U^{1,c'}\equiv F_V^{1,c'}$. Let $P_1=\{v_2v_3v_4v_5\colon v\in U^{1,c}\}$ and  $P_2=\{v_2v_3v_4v_5\colon v\in U^{1,c'}\}$. 
It is easy to check that we may pass from $P_1$ to $P_2$ for each form of $F_V^{1,c}\cup F_V^{1,c'}$. (At the end of the proof we show how to pass from $P_1$ to $P_2$ for the codes number one (Figure 5) and four.) 
Consequently, we can pass from $U^{1,c}\cup U^{1,c'}$ to $cP_2\cup U^{1,c'}$, where $cP_2=\{cv\colon v\in P_2\}$. The code $cP_2\cup U^{1,c'}$ consists of twin pairs of the form $cv,c'v$, where $v\in P_2$. Thus, we may pass by the glue and cut procedure from $cP_2\cup U^{1,c'}$ to $W$, where $W$ contains only twin pairs of the form $bv,b'v$ for $v\in P_2$. Clearly, $W\equiv F_V^{1,c}\cup F_V^{1,c'}$. Let $Q=F_V\setminus (F_V^{1,c}\cup F_V^{1,c'})\cup W$. Note that $Q^{1,c}\cup Q^{1,c}=\emptyset$. Thus, repeating such reductions we may pass from $V$ to $Q$ by the glue and cut procedure, where $Q\subset \{a,a',b,b',*\}^5$ and therefore, we may pass by the glue and cut procedure from $V$ to $\bar{Q}\subset \{a,a',b,b'\}^5$.


\vspace{-0mm}
{\center
\includegraphics[width=13cm]{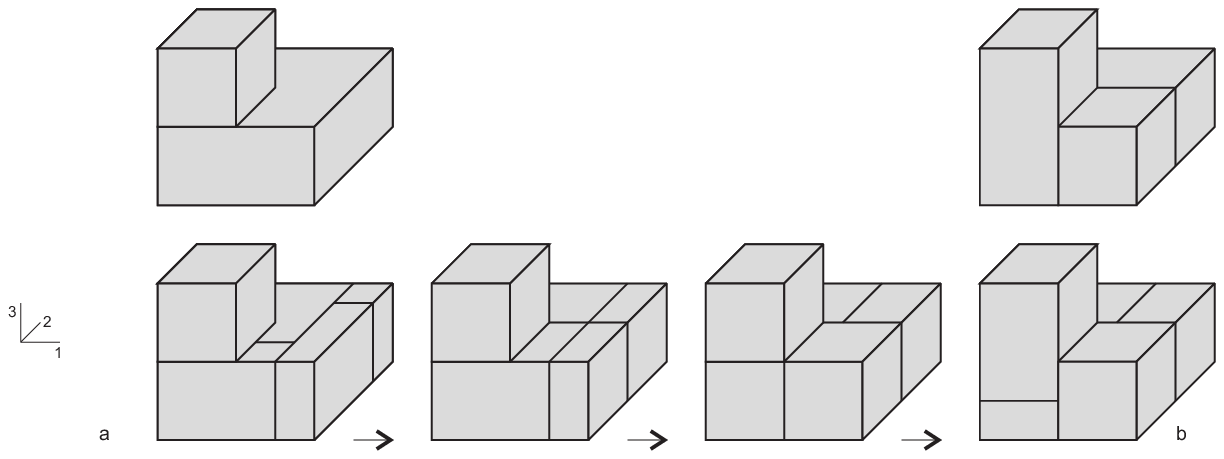}\\
}

\medskip
\noindent{\footnotesize Figure 5: Skipping in $F_V^{1,c}$ and $F_V^{1,c'}$ the first two letters we obtain the codes $A=\{**a, aaa'\}$, $B=\{aa*, a'*a, aa'a\}$, respectively. A realization of $A$  is given in the top left corner, and of $B$ in the top right corner. In picture ($a$) we have a realization of the code $U^{1,c}$ in which two first letters are skipped and similarly the picture ($b$) presents the code $U^{1,c'}$ with two first letters removed.      
}

\medskip
Observe that, the codes $P_1,P_2$ for the code $F_V^{1,c}\cup F_V^{1,c'}$ at position four cannot be represent in three dimensions. To  show how to pass from $P_1$ to $P_2$ let $p,q,s\in S$, and let $P_1=\{asa'b,as'a'b,qaaa,q'aaa,aaaa'\}$. We have $\{asa'b,as'a'b\}\rightarrow \{aaa'b,aa'a'b\}$ and  $\{qaaa,q'aaa\}\rightarrow \{aaaa,a'aaa\}$. Next we pass from $\{aaaa,aaaa'\}$ to $\{aaab,aaab'\}$. Thus, we have been passed by the glue and cut procedure from $P_1$ to the code $\{aaa'b, aa'a'b, aaab,aaab', a'aaa\}$. Now $\{aaa'b,aaab\}\rightarrow \{aap'b,aapb\}$, which gives $P_2=\{aap'b,aapb, aa'a'b, aaab', a'aaa\}$.
\hfill{$\square$}




  



\end{document}